\documentclass[a4paper,10pt]{article}
\usepackage{amsmath,amssymb}
\usepackage[latin1]{inputenc}
\usepackage{mathrsfs}

\newcommand{\CC}{\mathbb{C}}

\renewcommand{\SS}{\mathsf{S}}
\newcommand{\TT}{\mathsf{T}}

\newcommand{\HH}{\mathbf{H}}
\newcommand{\MM}{\mathbf{M}}

\newcommand{\ri}{\mathscr{R}}
\newcommand{\lf}{\mathscr{L}}

\newtheorem{theorem}{Theorem}[section] 
\newtheorem{proposition}[theorem]{Proposition} 
\newtheorem{conjecture}[theorem]{Conjecture} 
\newtheorem{corollary}[theorem]{Corollary} 
\newtheorem{lemma}[theorem]{Lemma}

\newenvironment{proof}{\begin{trivlist}\item{\bf{Proof.}}}
  {\hfill\rule{2mm}{2mm}\end{trivlist}}

\title{Cohomology rings of toric varieties \\ assigned to cluster
  quivers: \\ the case of unioriented quivers of type $A$} \author{F.
  Chapoton}

\date{\today}

\begin{document}

\maketitle

\begin{abstract}
  The theory of cluster algebras of S. Fomin and A. Zelevinsky has
  assigned a fan to each Dynkin diagram. Then A. Buan, R. Marsh, M.
  Reineke, I. Reiten and G. Todorov have generalized this construction
  using arbitrary quivers on Dynkin diagrams. In the special case of
  the unioriented quiver of type $A$, we describe the cohomology ring
  of the toric variety associated to this fan. A natural base is
  obtained and an explicit rule is given for the product of any two
  generators.
\end{abstract}

\setcounter{section}{-1}

\section{Introduction}

Cluster algebras were introduced by S. Fomin and A.  Zelevinsky
\cite{clu1,clu2,ysyst} for studying dual canonical bases in quantum
groups and total positivity in Lie groups. One important result of
this theory is the classification theorem of cluster algebras of
finite type by the Killing-Cartan list of root systems.  Part of the
proof consists of building a cluster algebra of finite type starting
from a given root system. In this construction, a simplicial fan is
associated with each finite root system, hence also a smooth toric
variety. It was proved in \cite{cfz} that these toric varieties are
projective.

Later, in the case of simply-laced Dynkin diagrams, this fan was seen
to be a special case of a construction starting from any quiver on the
given Dynkin diagram \cite{bmrrt}. The fan studied before corresponds in
this setting to the alternating quiver. It would be worth proving that
all the associated toric varieties of these quiver fans, known to be
smooth, are also projective.

In fact, this may be only the tip of something. There should be a
systematic way to define a fan starting from any cluster in a cluster
algebra of finite type, in such a way that the associated toric
variety is smooth and projective. This construction should of course
recover the preceding one, in case the chosen cluster is associated
with a quiver.

Moreover, based on experimental evidence, these fans should have the
following properties. First, they should not only be simplicial but
even smooth, meaning that each cone is spanned by an integral base.
Then the cone of ample divisors in the second integer cohomology group
of the toric variety should be smooth in the same sense. If true, this
would provide a natural base of this cohomology group and therefore a
natural set of generators of the cohomology ring. Then there should be
a quadratic presentation of the cohomology ring and a base of the
cohomology ring consisting of monomials in the distinguished
generators.

All these properties have been checked for low-dimensional alternating
quivers of type $A$. The aim of the present article is to prove part
of these statements in the case of unioriented quivers of type $A$.

More precisely, we obtain a base and a quadratic presentation of the
cohomology ring. The distinguished generators should be the extremal
vectors of the cone of ample divisors, but we do not prove that here.
To say the truth, this was however the way we guessed them by looking
at low-dimensional cases.

Let us remark that the rings studied here have some obvious similarity
with some rings related to the hyperplane arrangement of a root system
and to non-nesting partitions, which were considered in
\cite{heaviside}. We do not known what should be the meaning of this
resemblance.

\section{A toric variety associated to a Dynkin quiver}

Let us fix an integer $n$ once and for all and denote by $[n]$ the set
$\{1,2,\dots,n\}$.

For each quiver of Dynkin type, a fan has been defined on the set of
almost positive roots in \cite{bmrrt}. A similar construction is expected
to hold starting from any seed in any cluster algebra of finite type.

Let us recall the construction of \cite{bmrrt} in the case of the
unioriented quiver of type $A_n$. So let $Q_n$ be the quiver with $n$
vertices and arrows from $i$ to $i+1$ for $1 \leq i < n$.

By a simple instance of a famous Theorem of Gabriel, indecomposable
modules over $Q_n$ are in bijection with positive roots in the $A_n$
root system. Recall that these positive roots are indexed by the
intervals $[i,j]$ in the set $[n]$. Here, the indecomposable module
associated to $[i,j]$ is given by the space $\CC$ on each vertex $k$
between $i$ and $j$, the null space elsewhere and identity maps when
possible.

Let us introduce some notation. Let $\Phi_{>0}$ be the set of positive
roots, \textit{i.e.} the set of intervals in $[n]$. The intervals
$[i,i]$ are called simple roots and the set of simple roots is denoted
by $\Pi$. Let $\Phi_{>1}$ be the set of non-simple roots. Let
$\Phi_{\geq -1}$ be the disjoint union of $\Phi_{>0}$ with a copy of
$\Pi$ denoted by $-\Pi$. The elements of $\Phi_{\geq -1}$ are called
almost positive roots and the elements of $-\Pi$ are called negative
simple roots. In the sequel, we will denote by Greek letters the roots
\textit{i.e.} elements of $\Phi_{>0}$ and by Latin letters
(corresponding to elements of $[n]$) the simple roots or their
opposite.

One says that $i \in \alpha$ if $\alpha=[j,k]$ and $j \leq i\leq k$.

Let us say that two roots $\alpha=[i,j]$ and $\beta=[k,\ell]$ in
$\Phi_{>0}$ \textit{overlap} if one has $i\leq k\leq j \leq \ell$ or
$k \leq i \leq \ell \leq j$. Let us say that they \textit{overlap
  strictly} if they overlap and neither $\alpha \subseteq \beta$ nor
$\beta \subseteq \alpha$.

To define the fan on the set of vectors $\Phi_{\geq -1}$, one needs a
symmetric binary relation on the set $\Phi_{\geq -1}$, called the
compatibility relation. The general definition is given in term of
\rm{Ext}-groups in a triangulated category called the cluster category,
which is defined as a quotient of the derived category of modules on
the chosen quiver. We will just state the result for the quiver $Q_n$.

An element $-i$ of $-\Pi$ is compatible with any other element $-j$ of
$-\Pi$.

An element $-i$ of $-\Pi$ is compatible $\alpha \in \Phi_{>0}$ if and
only if $i\not\in\alpha$.

Two elements $\alpha$ and $\beta$ in $\Phi_{>0}$ are not compatible if
and only if
\begin{enumerate}
\item[(i)] either $\alpha \cap \beta = \emptyset$ and $\alpha \cup
  \beta \in \Phi_{>0}$ (adjacent roots),
\item[(ii)] or $\alpha \not \subseteq \beta$, $\beta \not \subseteq
  \alpha$ and $\alpha \cap \beta \not= \emptyset$ (strictly
  overlapping roots).
\end{enumerate}

Here comes the description of the fan $\Sigma(Q_n)$. First we map the
set $\Phi_{\geq -1}$ into the free abelian group generated by variables
$\{\alpha_1,\dots,\alpha_n\}$ by
\begin{equation}
  \begin{cases}
    -i \mapsto -\alpha_i,\\
    \alpha \mapsto \sum_{i \in \alpha} \alpha_i.
  \end{cases}
\end{equation}
then a subset of (the image of) $\Phi_{\geq -1}$ spans a cone of
$\Sigma(Q_n)$ if and only if it is made of pairwise compatible elements.
 
It is known that the number of maximal cones of $\Sigma(Q_n)$ is the
number of clusters of type $A_n$, which is the Catalan number
\begin{equation}
  c_{n+1}=\frac{1}{n+2}\binom{2n+2}{n+1}.
\end{equation}
 
Let us define $\max \alpha$ and $\min \alpha$ for $\alpha=[i,j] \in
\Phi_{>0}$ to be $i$ and $j$ respectively. Let us define $\ri \alpha$
and $\lf \alpha$ for $\alpha=[i,j] \in \Phi_{>1}$ to be $\alpha
\setminus \min \alpha$ and $\alpha \setminus \max \alpha$
respectively.

For $\ell \in [n]$ and $\alpha=[i,j] \in \Phi_{>1}$ with $\ell \in
\lf\alpha$, we define $\ell /\alpha$ to be the root $[\ell,j]$ in
$\Phi_{>1}$ obtained by cutting the left-hand side of $\alpha$.
Similarly, if $\ell \in \ri \alpha$, let $\alpha /\ell$ be the root
$[i,\ell]$ in $\Phi_{>1}$ defined by cutting the right-hand side of
$\alpha$.

\section{Standard presentation of the cohomology}

Let $\Sigma$ be a smooth complete fan. Then there exists a standard
description of the integer cohomology ring of the smooth toric variety
$X_\Sigma$, see \cite{danilov} and \cite[\S 5.2]{fulton}. Let us
recall briefly this construction.

The cohomology ring $\HH^*(X_\Sigma)$ is generated by variables
$\TT_u$ indexed by the set of $1$-dimensional cones in the fan
$\Sigma$. Then there are linear and quadratic relations between these
generators. The linear relations are 
\begin{equation}
  \sum_{u} \langle v,w_u \rangle \TT_u =0,
\end{equation}
where $v$ runs through a base of the dual lattice and $w_u$ is the
unique integral generating vector for the cone $u$. The quadratic
relations are the following: the product $\TT_u \TT_v$ vanishes as
soon as there is no cone $\sigma$ in $\Sigma$ containing both $u$ and
$v$.

It is also known that $\HH^*(X_\Sigma)$ is a free abelian group of
rank the number of maximal cones of $\Sigma$ \cite[Th. 10.8]{danilov}.

For the fans $\Sigma(Q_n)$ we are interested in, this amounts to the
following description. 

\begin{proposition}
  \label{presH}
  The cohomology ring $\HH^*(X_{\Sigma(Q_n)})$ is presented by the
  generators $\TT_{-i}$ for $i\in[n]$ and $\TT_\alpha$ for $\alpha \in
  \Phi_{>0}$, the linear relations
  \begin{equation}
    \label{linrel}
    \TT_{-i}=\sum_{i\in \alpha \in \Phi_{>0}} \TT_\alpha\quad \text {  for
    } \,i\in[n],
  \end{equation}
  and the quadratic relations
  \begin{equation}
    \label{quadHi}
    \TT_{-i} \TT_\alpha=0\quad \text {  when 
    }\quad i \in \alpha \in \Phi_{>0},
  \end{equation}
  and
  \begin{equation}
    \label{quadHalpha}
    \TT_\alpha \TT_\beta=0
  \end{equation}
  when
  \begin{enumerate}
  \item[(i)] either $\alpha \cap \beta = \emptyset$ and $\alpha \cup
    \beta \in \Phi_{>0}$ (adjacent roots),
  \item[(ii)] or $\alpha \not \subseteq \beta$, $\beta \not \subseteq
    \alpha$ and $\alpha \cap \beta \not= \emptyset$ (strictly
    overlapping roots).
  \end{enumerate}

  The rank of the free abelian group $\HH^*(X_{\Sigma(Q_n)})$ is the
  Catalan number $c_{n+1}$.
\end{proposition}

When $\alpha=[i,j]\in\Phi_{>1}$, we will sometimes denote $\TT_\alpha$
by $\TT_{i,j}$.

One can rewrite the quadratic relations involving the variables
$\TT_i$ by using the linear relations (\ref{linrel}) to eliminate
these variables.

The relation $\TT_i \TT_{i+1}=0$ becomes
\begin{equation}
  \label{longue_ii}
  \TT_{-i}\TT_{-i-1}-\sum_{i+1 < k'}\TT_{-i}\TT_{i+1,k'}
  -\sum_{j<i}\TT_{j,i}\TT_{-i-1} +\sum_{{{j \leq i \leq k ; j<k
      }\atop{j' \leq i+1 \leq k' ; j'< k'}}\atop{\text{inclusion}}} \TT_{j,k}\TT_{j',k'}=0,
\end{equation}
where ``inclusion'' means that either $[j,k]$ is contained in $[j',k']$ or
vice-versa.

The relation $\TT_i \TT_{\alpha}=0$, with $\alpha\in \Phi_{>0}$, $i
\not \in \alpha$ and $\alpha$ containing either $i+1$ or $i-1$, becomes
\begin{equation}
  \label{longue_bout}
  \TT_{-i} \TT_{\alpha}=\sum_{{i \in \beta}\, ;\, {\alpha \subseteq \beta}} \TT_{\alpha}\TT_{\beta}.
\end{equation}

\begin{lemma}
  \label{carre_nul}
  The square of $\TT_{-i}$ is zero for all $i \in [n]$.
\end{lemma}
\begin{proof}
  One has
  \begin{equation}
    \TT_{-i}^2=\TT_{-i} \left(\sum_{i\in \alpha \in \Phi_{>0}} \TT_\alpha\right), 
  \end{equation}
  which vanishes by relations (\ref{quadHi}).
\end{proof}

\section{A ring with a quadratic presentation}

\subsection{Presentation}

Let us introduce a ring $\MM^*(n)$. Our aim will be to show that
this ring is isomorphic to $\HH^*(X_{\Sigma(Q_n)})$.

The ring $\MM^*(n)$ is the commutative ring generated by
variables $\SS_i$ for $i\in[n]$ and $\SS_\alpha$ for $\alpha \in
\Phi_{>1}$, modulo the following relations:
\begin{equation}
  \label{rela1}
  \SS_i^2=0, \qquad \SS_i \SS_\alpha=\sum_{j \in \alpha, j\not= i} \SS_i \SS_j \quad
  \text{when  } \,i\in \alpha,
\end{equation}
 and 
\begin{equation}
  \label{rela2}
 \SS_\alpha \SS_\beta=\sum_{i < j \in \alpha \cap \beta} \SS_i\SS_j
 + \sum_{\ell\in \ri\alpha \cap \lf\beta} \SS_{\alpha/\ell}\SS_{\ell/\beta}-
 \sum_{\ell,\ell+1 \in \ri\alpha \cap \lf\beta}\SS_{\alpha/\ell}\SS_{\ell+1 /\beta},
\end{equation}
whenever $\alpha$ and $\beta$ overlap with $\alpha \cap \beta$ of
cardinal at least $2$ (one assumes that $\alpha$ is on the left of
$\beta$).

Remark: the ring $\MM^*(n)$ is obviously graded with generators
of degree one. 

When $\alpha=[i,j]\in\Phi_{>1}$, we will sometimes denote $\SS_\alpha$
by $\SS_{i,j}$.

\subsection{Combinatorial preliminaries: codes and $U$-sets}

A \textit{code} is a word $C$ of length $n$ in the alphabet
$\{L,R,L\!R,V\}$ such that
\begin{itemize}
\item It contains as many letters $L$ as letters $R$.
\item Any left prefix contains at least as many letters $L$ as letters
  $R$.
\end{itemize}

Note that $L$ is for ``links'', $R$ for ``recht'' and $V$ for ``vakuum''.

The \textit{degree} $\deg(C)$ of a code $C$ is the number of symbols
$L$ seen in the word, \textit{i.e.} the number of letters $L$ plus the
number of letters $L\!R$. There is a natural duality operation $C
\mapsto C^*$ on codes given by the replacement of all occurrences of
$L\!R$ by $V$ and vice-versa. This involution maps a code of degree
$k$ to a code of degree $n-k$. Hence there is a unique code of length
$n$ and degree $n$, made of $n$ letters $L\!R$.

It should be a simple combinatorial exercise for the reader to check
that the number of codes of length $n$ is the Catalan number $c_{n+1}$.

\medskip

A $U$-set is a subset $u$ of $[n] \sqcup \Phi_{>1}$ such that
\begin{itemize}
\item[(i)] If $i\in u$ and $\alpha \in u$, then $i\not\in \alpha$.
\item[(ii)] If $\alpha$ and $\beta$ in $u$ are overlapping, then
  $\alpha\cap \beta$ is a singleton.
\end{itemize}

Then $U$-sets are in bijection with codes as follows. A $U$-set $u$ is
mapped to the code $C$ obtained by writing a $L$ at position $i$ for
each non-simple root $\alpha$ starting at $i$ in $u$, a $R$ at
position $i$ for each non-simple root $\alpha$ ending at $i$ in $u$, a
$L\!R$ at position $i$ for each $i$ in $u$ and then filling the word
with $V$. Note that the letter $L\!R$ can either be obtained directly
as such or as the successive writing of $L$ and $R$ at the same place.

For example, the $U$-set $\{[1],[3,4],[4,6],[6,7]\}$ is mapped to the code 
\begin{equation}
  (L\!R)(V)(L)(L\!R)(V)(L\!R)(R).
\end{equation}

The reverse bijection from codes to $U$-sets is easy and left to the
reader. 

By this correspondence between codes and $U$-sets, the degree of a
code is mapped to the cardinality of the associated $U$-set. There is
an induced duality on $U$-sets which will be used later.

\subsection{Spanning set}

We want to show that there is a spanning set of $\MM^*(n)$ indexed
by $U$-sets. First for each $U$-set $u$, one can define a monomial
$\SS^u$ in $\MM^*(n)$ as the product of variables $\SS_i$ and
$\SS_\alpha$ over the elements of $u$.

Let the height of any monomial in $\MM^*(n)$ be the sum of the
height of its variables, where the generator $\SS_i$ has height $1$
and $\SS_\alpha$ has height $\#\alpha$.

\begin{lemma}
  \label{span}
  The ring $\MM^*(n)$ is spanned by the monomials $\SS^u$, where $u$
  runs over the set of $U$-sets.
\end{lemma}

\begin{proof}
  Using the defining relations (\ref{rela1}) and (\ref{rela2}) of
  $\MM^*(n)$, one can replace any monomial not of the form $\SS^u$ for
  some $U$-set $u$ by a linear combination of monomials of strictly
  smaller height. The Lemma follows by induction on height.
\end{proof}

\section{Isomorphism and consequences}

\subsection{Isomorphism}

Let us now describe a map $\psi$ from $\MM^*(n)$ to
$\HH^*(X_{\Sigma(Q_n)})$ and prove that it is an isomorphism.

Define $\psi$ on the generators of $\MM^*(n)$ by
\begin{equation}
  \label{psii}
  \psi(\SS_i)=\TT_{-i} \quad\text{for}\quad i\in[n],
\end{equation}
and
\begin{equation}
  \label{psialpha}
  \psi(\SS_\alpha)=\sum_{i \in \alpha} \TT_{-i}-\sum_{\alpha\subseteq
    \beta} \TT_{\beta}\quad\text{for}\quad \alpha\in\Phi_{>1}. 
\end{equation}

\begin{proposition}
  Formulas (\ref{psii}) and (\ref{psialpha}) define a morphism of
  rings $\psi$ from $\MM^*(n)$ to $\HH^*(X_{\Sigma(Q_n)})$. 
\end{proposition}

\begin{proof}
  Let us first check that relations (\ref{rela1}) hold. By
  (\ref{psii}), one has
  \begin{equation}
    \psi(S_i^2)=
    \TT_{-i}^2, 
  \end{equation}
  which vanishes by Lemma \ref{carre_nul}. One also has
  \begin{equation}
    \psi(\SS_i \SS_\alpha)=
    \TT_{-i} \left(\sum_{j \in \alpha} \TT_{-j}-\sum_{\alpha\subseteq
    \alpha '} \TT_{\alpha'}\right). 
  \end{equation}
  By relations (\ref{quadHi}), this becomes, as expected,
  \begin{equation}
    \TT_{-i} \sum_{j \in \alpha, j \not=i} \TT_{-j}
    =\psi\left(\sum_{j \in \alpha, j \not=i} \SS_i \SS_{j}\right). 
  \end{equation}

  \smallskip

  Let us now check that relations (\ref{rela2}) hold. It is necessary
  to distinguish two cases.

  First consider the case when $\ri \alpha \cap \lf \beta$ is empty.
  One can show that this implies that $\alpha$ and $\beta$ are the
  same $[i,i+1]$ for some $i$. One has to check the vanishing of the
  image by $\psi$ of
  \begin{equation}
    \SS_{i,i+1}^2-\SS_{i}\SS_{i+1}.
  \end{equation}
  This is given by
  \begin{equation}
    \left(\TT_{-i}+\TT_{-i-1}-\sum_{j \leq i < i+1 \leq k} \TT_{j,k}\right)^2-\TT_{-i}\TT_{-i-1}.
  \end{equation}
  By relations (\ref{quadHi}), this is
  \begin{equation}
    \TT_{-i}\TT_{-i-1}+\left(\sum_{j \leq i < i+1 \leq k} \TT_{j,k}\right)^2.
  \end{equation}
  Using relations (\ref{quadHalpha}), this becomes
  \begin{equation}
    \TT_{-i}\TT_{-i-1}+\sum_{{{j \leq i < i+1 \leq k}\atop{j' \leq i <
        i+1 \leq k'}}\atop{\text{inclusion}}} \TT_{j,k}\TT_{j',k'} .
  \end{equation}
  where ``inclusion'' means that either $[j,k]$ is contained in $[j',k']$
  or vice-versa. Then using relation (\ref{longue_ii}) to eliminate
  $\TT_{-i}\TT_{-i-1}$, one gets
  \begin{equation}
    \sum_{i+1 < k'}\TT_{-i}\TT_{i+1,k'}
  +\sum_{j<i}\TT_{j,i}\TT_{-i-1}
  -\sum_{{j' \leq i \leq i+1 \leq k'}\atop{j' \leq j < i \leq k'}}\TT_{j,i}\TT_{j',k'}-\sum_{{j \leq i \leq k}\atop{j \leq i+1 < k' \leq k}}\TT_{i+1,k'}\TT_{j,k}.
  \end{equation}
  Then using relations (\ref{longue_bout}), the first and fourth term
  annihilate as do the second and third term.

  \medskip
  
  Let us now consider the case when $\ri \alpha \cap \lf \beta$ is
  not empty. We have to prove the vanishing of the image by $\psi$ of
  \begin{equation}
    \SS_\alpha \SS_\beta-\sum_{i < j \in \alpha \cap \beta} \SS_i\SS_j
 -\sum_{\ell\in \ri\alpha \cap \lf\beta} \SS_{\alpha/\ell}\SS_{\ell/\beta}+
 \sum_{\ell,\ell+1 \in \ri\alpha \cap \lf\beta}\SS_{\alpha/\ell}\SS_{\ell+1 /\beta}.
  \end{equation}

  This is
  \begin{multline}
    \label{grosse}
    \left(\sum_{i \in \alpha} \TT_{-i}-\sum_{\alpha\subseteq
        \alpha '} \TT_{\alpha'}\right)
    \left(\sum_{j \in \beta} \TT_{-j}-\sum_{\beta\subseteq
        \beta '} \TT_{\beta'}\right)
-\sum_{i<j \in \alpha \cap \beta} \TT_{-i}\TT_{-j}
\\    - \sum_{\ell \in \ri\alpha \cap \lf\beta} \left(\sum_{{i \in
          \alpha}\atop{i \leq \ell}} \TT_{-i}-\sum_{\alpha/\ell\subseteq
    \alpha '} \TT_{\alpha'}\right)\left(\sum_{{j \in \beta}\atop{j
      \geq \ell}}\TT_{-j}-\sum_{\ell/\beta\subseteq
    \beta '} \TT_{\beta'}\right)\\+
\sum_{\ell,\ell+1 \in \ri\alpha \cap \lf\beta} \left(\sum_{{i \in \alpha}\atop{i \leq \ell}} \TT_{-i}-\sum_{\alpha/\ell\subseteq
    \alpha '} \TT_{\alpha'}\right)\left(\sum_{{j \in \beta}\atop{j
      \geq \ell+1}}\TT_{-j}-\sum_{\ell+1/\beta\subseteq
    \beta '} \TT_{\beta'}\right).
  \end{multline}
  
  In this sum, consider first the terms of the shape
  $\TT_{-\star}\TT_{-*}$. Let us prove that their sum vanishes. First,
  using the fact that $\alpha$ and $\beta$ overlap with $\alpha$ on
  the left, and reversing summations, one gets
  \begin{equation}
    \label{abcde_gauche}
    \sum_{{{i \in
          \alpha}\atop{j \in \beta}}\atop{i<j}}\TT_{-i}\TT_{-j}- \sum_{{{i \in
          \alpha}\atop{j \in \beta}}\atop{i<j}}\sum_{\ell \in \ri\alpha
        \cap \lf\beta \cap [i,j]} \TT_{-i}\TT_{-j} +\sum_{{{i \in
          \alpha}\atop{j \in \beta}}\atop{i<j}}\sum_{\ell,\ell+1 \in \ri\alpha
        \cap \lf\beta \cap [i,j]} \TT_{-i}\TT_{-j}.
  \end{equation}
  Then it is enough to show that $\ri\alpha \cap \lf\beta \cap [i,j]$
  is not empty. This is clear if $i+1<j$, for any $i<k<j$ will do the
  job. Then if $j=i+1$, the intersection can be empty only if
  $\ri\alpha \cap \lf\beta $ is already empty, which is excluded by
  hypothesis.

  \medskip
 
  Then consider the terms of the shape $\TT_{-\star}\TT_{\alpha'}$ in
  (\ref{grosse}). Using the left-right symmetry of the situation, let
  us compute only the terms of the shape $\TT_{-j}\TT_{\alpha'}$ where
  $\alpha \subseteq \alpha'$ and $j \in \beta$. After reversal of
  summations, this sum is
  \begin{equation}
    -\sum_{{{j \in
      \beta}\atop{\alpha\subseteq
        \alpha '}}\atop{j\not\in\alpha'}} \TT_{-j}\TT_{\alpha'}+
    \sum_{{{j \in
      \beta}\atop{\min\alpha\in
        \alpha '}}\atop{j\not\in\alpha'}} \sum_{{\ell \in \ri\alpha
        \cap \lf\beta}\atop{\ell\in\alpha'}}\TT_{-j}\TT_{\alpha'}
-    \sum_{{{j \in
      \beta}\atop{\min\alpha\in
        \alpha '}}\atop{j\not\in\alpha'}} \sum_{{\ell,\ell+1 \in \ri\alpha
        \cap \lf\beta}\atop{\ell\in\alpha'}}\TT_{-j}\TT_{\alpha'}.
  \end{equation}
  The sum of the last two terms under the additional assumption that
  $\alpha \subseteq \alpha'$ annihilates with the first term. Here we
  used that $\ri\alpha \cap \lf\beta $ is not empty. Let us therefore
  assume that $\alpha \not\subseteq \alpha'$ in the two right terms.
  This means that $\max \alpha \not \in \alpha'$. Then both terms
  vanish unless $\alpha'$ meets $\ri\alpha \cap \lf\beta $. In this
  case, the sum vanishes unless $\max \alpha'=\max \ri\alpha \cap
  \lf\beta $. 

  This situation is possible if and only if $\max \alpha=\max \beta$,
  in which case one gets
  \begin{equation}
    \label{sub1}
    \sum_{\alpha \cap \alpha'=\lf \alpha} \TT_{-\max \alpha} \TT_{\alpha'}.
  \end{equation}
  A similar proof for the left-right symmetric summation, gives that
  the corresponding sum vanishes unless $\min \alpha=\min \beta$, in
  which case it is given by
  \begin{equation}
    \label{sub2}
    \sum_{\beta \cap \beta'=\ri \beta} \TT_{-\min \beta} \TT_{\beta'}.
  \end{equation}

  \medskip
   
  Then, at last, consider the terms of the shape
  $\TT_{\alpha'}\TT_{\beta'}$ in (\ref{grosse}). This is given by
  \begin{equation}
    \sum_{{{\alpha\subseteq
      \alpha '}\atop {\beta\subseteq
      \beta '}}\atop{\text{inclusion} }} \TT_{\alpha'} \TT_{\beta'}-
\sum_{{{\min \alpha \in \alpha'}\atop{\max\beta \in
      \beta'}}\atop{\text{inclusion}}}\sum_{{\ell \in \ri \alpha \cap
    \lf \beta}\atop{\ell \in \alpha'\cap\beta'}}\TT_{\alpha'}
\TT_{\beta'}+
\sum_{{{\min \alpha \in \alpha'}\atop{\max\beta \in
      \beta'}}\atop{\text{inclusion}}}\sum_{{{\ell,\ell+1 \in \ri \alpha \cap
    \lf \beta}\atop{\ell \in \alpha'}}\atop{\ell+1\cap\beta'}}\TT_{\alpha'}
\TT_{\beta'}.
  \end{equation}
  In each term, as $\alpha'\cap \beta'$ is necessarily not empty, the
  summation on $\alpha'$ and $\beta'$ can be restricted using
  relations (\ref{quadHalpha}) to the cases where $\alpha'\subseteq
  \beta'$ or vice-versa. This is the meaning of the ``inclusion''
  subscripts. The sum of the last two terms under the additional
  assumption that $\alpha \subseteq \alpha'$ and $\beta \subseteq
  \beta'$ annihilates with the first term. Here we used once again
  that $\ri\alpha \cap \lf\beta $ is not empty.

  Then one can assume in the right two terms that either $\alpha
  \not\subseteq \alpha'$ or $\beta \not\subseteq \beta'$. It turns
  out that these possibilities exclude each other because of the
  inclusion $\alpha'\subseteq \beta'$ or vice-versa.

  Let us compute the sum when $\alpha \not\subseteq \alpha'$, $\alpha
  \cup \beta \subseteq \beta'$ and $\alpha'\subseteq \beta'$. This is
  given by
  \begin{equation}
    -\sum_{{{\min \alpha \in \alpha'}\atop{\max\beta \in
      \beta'}}\atop{\alpha'\subseteq \beta'}}\sum_{{\ell \in \ri \alpha \cap
    \lf \beta}\atop{\ell \in \alpha'}}\TT_{\alpha'}\TT_{\beta'}+
\sum_{{{\min \alpha \in \alpha'}\atop{\max\beta \in
      \beta'}}\atop{\alpha'\subseteq \beta'}}\sum_{{\ell,\ell+1 \in \ri \alpha \cap
    \lf \beta}\atop{\ell \in \alpha'}}\TT_{\alpha'}\TT_{\beta'}.
  \end{equation}
 Then both terms
  vanish unless $\alpha'$ meets $\ri\alpha \cap \lf\beta $. In this
  case, the sum vanishes unless $\max \alpha'=\max \ri\alpha \cap
  \lf\beta $. 

  This situation is possible if and only if $\max \alpha=\max \beta$,
  in which case one gets
  \begin{equation}
    \label{sub3}
    \sum_{\alpha \cap \alpha'=\lf \alpha}\sum_{\alpha'\cup \alpha \subseteq \beta'} \TT_{\beta'} \TT_{\alpha'}.
  \end{equation}

  Similarly, the sum  when $\alpha \cup \beta \subseteq \alpha'$,
  $\beta \not \subseteq \beta'$ and $\beta'\subseteq \alpha'$ vanish
  unless $\min \alpha=\min \beta $, in which case it is given by 
  \begin{equation}
    \label{sub4}
    \sum_{\beta \cap \beta'=\ri \beta}\sum_{\beta'\cup \beta \subseteq \alpha'} \TT_{\alpha'} \TT_{\beta'}.
  \end{equation}
  
  Then gathering the terms (\ref{sub1}),(\ref{sub3}) and the terms
  (\ref{sub2}),(\ref{sub4}) and using relations (\ref{longue_bout}),
  one gets the expected vanishing of (\ref{grosse}) in all cases.
\end{proof}

\begin{theorem}
  The morphism $\psi$ from $\MM^*(n)$ to $\HH^*(X_{\Sigma(Q_n)})$ is
  an isomorphism.
\end{theorem}

\begin{proof}
  Let us first prove that $\psi$ is surjective. First it is clear from
  (\ref{psii}) that each $\TT_{-i}$ is in the image of $\psi$. Then
  one can see by Möbius inversion on (\ref{psialpha}) that each
  $\TT_{\alpha}$ for $\alpha\in\Phi_{>1}$ is also in the image of
  $\psi$. But these variables together generates
  $\HH^*(X_{\Sigma(Q_n)})$ because of the linear relations
  (\ref{linrel}).

  Now the ring $\HH^*(X_{\Sigma(Q_n)})$ is a free abelian group of rank
  the Catalan number $c_{n+1}$. By Lemma \ref{span}, the surjectivity
  of $\psi$ then implies that the monomials $S^u$, for $u$ in the set
  of $U$-sets, are linearly independent in $\MM^*(n)$. Hence they
  form a base of $\MM^*(n)$ and their images must be a base of
  $\HH^*(X_{\Sigma(Q_n)})$. So $\psi$ is an isomorphism.

\end{proof}

\subsection{Consequences}

The first consequence of this isomorphism is of course that the
monomials $S^u$ for $U$-sets $u$ form a base of $\MM^*(n)$. Let us
call it the natural base. From now on, we will identify $\MM^*(n)$
with $\HH^*(X_{\Sigma(Q_n)})$ by the mean of $\psi$.

There is a unique element of degree $n$ in the natural base, which is
the product of all $\SS_i$.

From Poincaré duality in the cohomology ring, one gets
\begin{corollary}
  The ring $\MM^*(n)$ is a graded Frobenius ring. 
\end{corollary}

\begin{theorem}
  The set of relations (\ref{rela1}) and (\ref{rela2}) is a
  (quadratic) Gröbner basis for the term order where variables of
  greater height are dominant.
\end{theorem}

\begin{proof}
  If this is not true, then there would exist another element in the
  Gröbner basis with a leading monomial of the form $\SS^u$ for some
  $U$-set $u$. This would contradict the fact that the monomials
  associated to $U$-sets are linearly independent.
\end{proof}

\begin{theorem}
  The ring $\MM^*(n)$ is Koszul as an associative algebra.
\end{theorem}

\begin{proof}
  This follows from the fact that it admits a quadratic Gröbner basis,
  see for example \cite{anick}.
\end{proof}

\begin{proposition}
  The ring $\MM^*(n)$ is filtered by the subspaces spanned by
  monomials $\SS^u$ of height less than a fixed bound.
\end{proposition}

\begin{proof}
  Indeed, the procedure of rewriting the product of two monomials in
  the natural base as a sum of elements of this base uses the
  Gröbner basis reduction, which can only decrease the height.
\end{proof}

\subsection{Duality between bottom $\TT$ and top $\SS$}

As said before, the natural base of $\MM^*(n)$ contains a unique
element of degree $n$, which is simply
\begin{equation}
  \prod_{i  \in [n]}\SS_{i}.
\end{equation}

The symmetric bilinear form $\langle\,,\, \rangle $ defining the
Frobenius structure of the graded ring $\MM^*(n)$ is given by the
coefficient of this unique element of degree $n$ in the expression in
the natural base of the product of two elements of $\MM^*(n)$.

By the graded Frobenius property, this bilinear map restricts to a
non-degenerate pairing between the subspace of degree $1$ (spanned by
generators) and the subspace of degree $n-1$.

Let us consider the natural base in degree $n-1$. It is indexed by
$U$-sets of cardinality $n-1$. Using the duality on $U$-sets coming
from the duality on codes, one can instead index this base by
$\Phi_{>0}$. Let $\SS'_{\alpha}$ be the element of the natural base in
degree $n-1$ assigned in this way to $\alpha\in \Phi_{>0}$.

\smallskip

By the Frobenius pairing, the natural base in degree $n-1$ has a
simple dual base in degree $1$:

\begin{proposition}
  The base $(\TT_\alpha)_{\alpha \in \Phi_{>0}}$ in degree $1$ is dual
  to the base $(\SS'_\alpha)_{\alpha \in \Phi_{>0}}$ in degree $n-1$
  for the Frobenius pairing: for all $\alpha,\beta$ in $\Phi_{>0}$, one has
  \begin{equation}
    \SS_\alpha=\sum_\beta \langle \SS_\alpha,\SS'_\beta\rangle \TT_\beta.
  \end{equation}
  
\end{proposition}

\begin{proof}
  The proof is based on the comparison between the explicit
  computation of the coefficients $\langle
  \SS_\alpha,\SS'_\beta\rangle$ and the change of base between the
  natural base $\SS_\alpha$ in degree $1$ and the base $\TT_\alpha$.

  Let us start with the change of basis between $\SS$ and $\TT$. Using
  (\ref{psii}), (\ref{psialpha}) and (\ref{linrel}), one finds that
  \begin{equation}
    \SS_i =\sum_{i \in \alpha \in \Phi_{>0}} \TT_\alpha,
  \end{equation}
  and for $\alpha\in \Phi_{>1}$, 
  \begin{equation}
    \SS_{\alpha}=\sum_{\alpha \subseteq
      \beta}(\#\alpha-1)\TT_\beta+\sum_{\alpha \not\subseteq \beta \in
      \Phi_{>0}}(\#\alpha \cap \beta)\TT_\beta.
  \end{equation}  
  Then it only remains to show that these formulas coincide with the
  value of the pairing. This is done below.
\end{proof}

Let us first state two useful Lemmas.
\begin{lemma}
  \label{minil1}
  For $1\leq i < j \leq n$, one has
  \begin{equation}
    \SS_{i,j} \left(\SS_{i,i+1} \dots \SS_{j-1,j}\right)=(j-i)\SS_i \dots \SS_j.
  \end{equation}
\end{lemma}
\begin{proof}
  This is a simple inductive computation in $\MM^*(n)$. This is easy
  if $j=i+1$. The inductive step first computes the product $\SS_{i,j}
  \SS_{i,i+1}$.
\end{proof}

\begin{lemma}
  \label{minil2}
  For $1\leq i < j \leq n$, one has
  \begin{equation}
    \SS_{i,j} \left(\SS_{i} \dots \SS_{j}\right)=0.
  \end{equation}
\end{lemma}
\begin{proof}
  Quite obvious from the defining relations, by induction.
\end{proof}

Let us now compute the pairing between elements of degree $1$ and
elements of degree $n-1$ in the natural base of $\MM^*(n)$.

\begin{proposition}
  The following equations hold for $\alpha,\beta$ in $\Phi_{>1}$:
  \begin{equation}
  \langle \SS_j,\SS'_i\rangle =\delta_{i=j},
\end{equation}
\begin{equation}
  \langle \SS_\beta,\SS'_i\rangle =\delta_{i\in\beta},
\end{equation}
\begin{equation}
  \langle \SS_j,\SS'_\alpha \rangle =\delta_{j \in \alpha},
\end{equation}
\begin{equation}
  \langle \SS_\beta,\SS'_\alpha \rangle = 
  \begin{cases}
    \# \beta -1 \text { if } \beta \subseteq \alpha,\\
    \# \alpha \cap \beta \text{ else.}
  \end{cases}
\end{equation}

\end{proposition}

\begin{proof}
  First, note that
  \begin{equation}
    \SS'_i=\prod_{j\not =i} \SS_j.
  \end{equation}
  This implies the first relation using that $\SS_i^2=0$ and the
  second relation using Lemma \ref{minil2}. Then note that for
  $\alpha=[i,j]$ with $i<j$, one has
  \begin{equation}
    \SS'_{\alpha}=(\SS_{i,i+1}\dots\SS_{j-1,j})\prod_{k\not \in \alpha}\SS_k.
  \end{equation}
  This easily implies the third relation. The fourth relation can be
  checked by distinguishing whether $\beta \subseteq \alpha$ or not
  and using Lemma \ref{minil1}.
\end{proof}

\section{Parabolic inclusions}

It follows from the presentation of the rings $\MM^*({n})$ that, for
any $n_1$ and $n_2$, there are morphisms of rings
\begin{equation}
  \MM^*({n_1})\otimes \MM^*({n_2}) \to \MM^*({n_1+n_2}),
\end{equation}
mapping the generators $\SS\otimes 1$ and $1\otimes \SS$ to some
generators $\SS$ according to the decomposition of the interval
$[n_1+n_2]$ into two consecutive intervals $[n_1]$ and $[n_2]$. 

These morphisms map the tensor product of the natural bases into the
natural base, hence they are injective. As the sum of the ranks is
smaller than the rank in general, they are not surjective. 

One can even see that, for a fixed $n$, the span of all the images of
these maps for varying $n_1,n_2$ of sum $n$ can not be the full ring
$\MM^*(n)$, for it can not contain the element $\SS_{1,n}$. One can
compute that the number of elements of the natural base which cannot
be reached in this way is the Catalan number $c_{n-1}$.

Through the isomorphisms with cohomology rings, these morphisms should
come from refinements of fans, inducing maps of toric varieties, hence
maps at the level of cohomology.

\section{Conjectural deformation}

It seems that one can replace the relations $\SS_i^2=0$ in the
presentation of $\MM^*(n)$ by the relation $\SS_i^2=\SS_i$ without
much harm. 

Let $\MM^{\text{def}}(n)$ be the commutative ring generated by
variables $\SS_i$ for $i\in[n]$ and $\SS_\alpha$ for $\alpha \in
\Phi_{>1}$, modulo the right half of relations (\ref{rela1}), all
relations (\ref{rela2}) and relations $\SS_i^2=\SS_i$.

Of course, this ring is not graded as one relation is no longer homogeneous.

\begin{conjecture}
  The ring $\MM^{\text{def}}(n)$ has dimension $c_{n+1}$.
\end{conjecture}

This has been checked by computer for $n \leq 6$. A strategy of proof
would be to show that this set of relations is still a Gröbner basis.
As part of this check, it is easy to see that the reduction of the
monomials $\SS_i^2 \SS_\alpha$ for $i \in \alpha \in \Phi_{>1}$ works well.

\bibliographystyle{plain}
\bibliography{cluster_coho}

\begin{thebibliography}{1}

\bibitem{anick}
David~J. Anick.
\newblock On the homology of associative algebras.
\newblock {\em Trans. Amer. Math. Soc.}, 296(2):641--659, 1986.

\bibitem{bmrrt}
Aslak~Bakke Buan, Robert Marsh, Markus Reineke, Idun Reiten, and Gordana
  Todorov.
\newblock {Tilting theory and cluster combinatorics}.
\newblock arXiv:math.RT/0402054.

\bibitem{heaviside}
Fr{\'e}d{\'e}ric Chapoton.
\newblock {Antichains of positive roots and Heaviside functions}.
\newblock arXiv:math.CO/0303220.

\bibitem{cfz}
Fr{\'e}d{\'e}ric Chapoton, Sergey Fomin, and Andrei Zelevinsky.
\newblock Polytopal realizations of generalized associahedra.
\newblock {\em Canad. Math. Bull.}, 45(4):537--566, 2002.
\newblock Dedicated to Robert V.\ Moody.

\bibitem{danilov}
V.~I. Danilov.
\newblock The geometry of toric varieties.
\newblock {\em Uspekhi Mat. Nauk}, 33(2(200)):85--134, 247, 1978.

\bibitem{clu1}
Sergey Fomin and Andrei Zelevinsky.
\newblock Cluster algebras. {I}. {F}oundations.
\newblock {\em J. Amer. Math. Soc.}, 15(2):497--529 (electronic), 2002.

\bibitem{clu2}
Sergey Fomin and Andrei Zelevinsky.
\newblock Cluster algebras. {II}. {F}inite type classification.
\newblock {\em Invent. Math.}, 154(1):63--121, 2003.

\bibitem{ysyst}
Sergey Fomin and Andrei Zelevinsky.
\newblock {$Y$}-systems and generalized associahedra.
\newblock {\em Ann. of Math. (2)}, 158(3):977--1018, 2003.

\bibitem{fulton}
William Fulton.
\newblock {\em Introduction to toric varieties}, volume 131 of {\em Annals of
  Mathematics Studies}.
\newblock Princeton University Press, Princeton, NJ, 1993.
\newblock The William H. Roever Lectures in Geometry.

\end{thebibliography}

\end{document}